\documentclass[14pt]{article}
\pdfoutput=1
\usepackage{authblk}
\usepackage[12pt]{extsizes}
\usepackage{geometry}
\usepackage[ruled,vlined]{algorithm2e}
\usepackage{amssymb}
\usepackage{graphicx}
\usepackage{float}

\usepackage{amsmath}
\usepackage{verbatim}

\def\keywords{\vspace{.5em}
{\textit{Keywords}:\,\relax%
}}

\title{
  High Order Hermite Finite Difference Method for Euler/Navier-Stokes Equations in 2D Unstructured Meshes
}

\author[a,b]{Zeyuan Zhou}
\author[a]{Mei-Yuan Zhen}
\author[a]{Kun Qu \thanks{corresponding author: kunqu@nwpu.edu.cn, School of Aeronautics, Northwestern Polytechnical University, Shaanxi Province, Xi'an, 710072, China}}
\author[a]{Jin-Sheng Cai}

\affil[a] {School of Aeronautics, Northwestern Polytechnical University}
\affil[b] {School of Mechanics, Northwestern Polytechnical University}

\date{\today}

\begin{document}
\maketitle


\begin{abstract}
  A high order finite difference method is proposed for unstructured meshes to simulate compressible inviscid/viscous flows with/without discontinuities. In this method, based on the strong form equation, the divergence of the flux on each vertex is computed directly from fluxes nearby by means of high order least-square. In order to capture discontinuities, numerical flux of high order accuracy is calculated on each edge and serves as supporting data of the least-square computation of the divergence. The high accuracy of the numerical flux depends on the high order WENO interpolation on each edge. To reduce the computing cost and complexity, a curvlinear stencil is assembled for each edge so that the economical one-dimensional WENO interpolation can be applied. With the derivatives introduced, two-dimensional Hermite interpolation on a curvilinear stencil is applied to keep the stencil compact and avoids using many supporting points. In smooth region, the Hermite least-square 2D interpolation of 5 nodes is adopted directly to achieve the fifth order accuracy. Near a discontinuity, three values obtained by means of least-square 2D interpolation of 3 nodes, are weighted to obtain one value of the second order accuracy. After obtaining the flow states on both sides of the middle point of an edge, numerical flux of high order accuracy along the edge can be calculated. For inviscid flux, analytical flux on vertices and numerical flux along edges are used to compute the divergence. While for viscous flux, only analytical viscous flux on vertices are used. The divergence of the fluxes and their derivatives on each vertex are used to update the conservative variables and their derivatives with an explicit Runger-Kutta time scheme. Several canonical numerical cases were solved to test the accuracy and the capability of shock capturing of this method.

  \keywords{finite difference method, high order method, unstructured mesh, WENO interpolation, least square}

  \end{abstract}

  \section{Introduction} \label{sec_intro}

  For more than twenty years, many researchers have been working on high order methods for computational fluid dynamics(CFD) since they are crucial in many important scientific and engineering fields, such as flow instabilty, turbulent flows, aeroacoustics, chemical reactive flows, multipahse flows ...
  So far, high order finite difference method(FDM), finite volume method(FVM), discontinuous Galerkin(DG), correction procedure via reconstruction(CPR), spectral difference(SD) and flux reconstruction(FR) are most popular high order methods in CFD community.
  But they are still far from mature. And any one of them has some disadvantages and problems to be fixed.

  For high order FVM, it relies on high order reconstruction on stencils of cells and integaration of flux on faces.
  With multidimensional WENO reconstruction, high order FVM can capture discontinuities effectively.
  But in this case, WENO reconstruction requires assembling multiple sub-stencils, which is very expensive and complicated, even not robust.
  At the same time, each cell of FVM only contains one DOF, which makes the cost per DOF very high and increases \emph{exponentially} with the space dimension.

  Although FR\cite{huynh_flux_2007} (including SD\cite{liu2006spectral,wang2007spectral}, CPR and DG \cite{cockburn1989tvb,cockburn2001runge} since they can be written in the form of FR) in unstructred meshes is more economical, they can not capture discontinuities because the shape functions in each FR element are smooth polynomials of high degrees.
  Although some techniques, such as artifical viscosity, slope limiter, sub-cells and shock-fitting coupled with mesh adaptation, are applied into FR to capture discontinuities, they make FR more complicated and expensive. 

  %
  %
  %
  %
  %
  %

  On the contrary, high order FDM is of very high efficiency and very low complexity because the computation is performed dimension by dimension.
  This makes the computing computing cost per DOF increase \emph{linearly} with the space dimension. 
  Nowadays, high order WENO finite difference scheme \cite{shu_efficient_1988,jiang_efficient_1996} and its branches are highly developed and applied widely in simulations of shock waves, vortices, shear layers, interfaces ...
  Of course, the disadvantage of high order FDM is obvious: it can only be applied in smooth structured grids which are difficult to generate for complex geometries so far.

  %
  %
  %
  %

  There is a question raised naturally: is that possible to apply the mature shock-capturing techniques developed in FDM into the contex of unstructured meshes?
  There is no successful story so far.
  The advantage of FDM is that there is special one-dimensional topology in multidimensional structured grids.
  But for FEM, FVM, DG and FR in unstructured meshes, multidimensional interpolation or reconstructuion was adopted since there is not such special one-dimensional topology.
  All these stories make people belive that applying multidimensional approximation in unstructured meshes is neccesary.

  \textbf{Can we find one-dimensional topological structures in multidimensional unstructured meshes?}
  If it is possible, we might be able to apply the mature shock-capturing techniques developed in FDM into the contex of unstructured meshes.
  In fact, we \textbf{can} find one-dimensional topological structures in multidimensional unstructured meshes at the first glance.
  \begin{figure}\label{fig:assemble-stencil}
    \centering
    \includegraphics[width=0.4\textwidth]{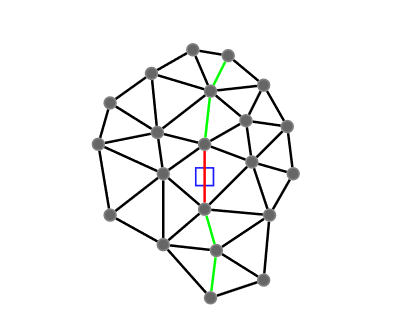}
    \caption{The red and green edges make of an one-dimensional structure in an unstructured mesh. The square symbol is the middle point of the green edge.}
  \end{figure}
  The connected edges in an unstructured mesh can be regarded as a curve which is an one-dimensional topological structure.
  It can be easily assmbled by selecting connected edges with smallest angle of deflection.
  From the six nodes on the five edges in Fig. \ref{fig:assemble-stencil}, we maybe able to interpolate variables on the middle point of the middle edge.
  Of course, it is not so smooth as the curves in a good structured grid, which may make error larger when applying FD schemes along it.
  But it is still a good begining to find a way to apply the mature shock-capturing techniques developed in FDM into the contex of unstructured meshes.

  In this work, the authors try to extend the discontinuity-capturing techniques in FDM to unstructured meshes in order to develop an efficient high order discontinuity-capturing method in unstructured meshes.
  In the next section, the one-dimensional WENO interpolation is directly extended to the curves connected by edges in unstructured meshes.
  But numerical tests show that the convergence rate of this interpolation is only the first order even if the high order one-dimensional WENO interpolation is used.
  Further analysis shows that in the context of refining unstructured meshes, interpolation along a curve based on one-dimensional high degree polynomials always degenerates to the first order convergence rate.
  This conclusion means that the multidimensional interpolation is neccsary.
  In order to make the interpolation as economical as possible, we try to keep the curilinear stencil while applying multidimensional interpolation on it.
  %
  Since increasing the number of supporting nodes is not feasible, we choose to use Hermite interpolation which rely on the function values and their gradients on each supporting node.
  Coupled with the weights of WENO, both smooth region and discontinuities can be handeled.
  Numerical tests verified the accuracy  of this special Hermite interpolation.
  In this way, the numerical flux on the middle point of each edge can be calculated.
  With the flux data on both vertices and edges, a multidimensional high order least-square based difference scheme is used to computed the divergence of fluxes and their gradients on each vertex.
  An explicit  Runge-Ketta time scheme is adopted to update the conservative variables and their gradients on each vertex.
  In the 3rd section, the accuracy order is tested with the case of inviscid isentropic vortex and the viscous Couete flow.
  The circular Sod problem is used to test the capability of discontinuity-capturing of this meshod.

  \section{Methodology} \label{sec_methd}
  \subsection{Extending One-Dimensional WENO Interpolation in Multidimensional Unstructured Meshes} \label{subsec_MDweno}
  Assembling a curvilinear strencil in unstructred meshes is not difficult.
  For an end point of an edge, we just serach all the other edges linked at the end point to find one with smallest deflection.
  Performing twice searches to extend the edge from each of its end points, we obtain a five-edge stencil in Fig. \ref{fig:assemble-stencil}.

  Our work is inspire by Weighted Compact Nonlinear Scheme (WCNS) of Deng\cite{deng_developing_2000}.
  In WCNS-E5 scheme, the WENO5 interpolation is used to obtain left and right flow states on each interface where numerical flux can be computed by some Riemann flux solver.
  And then the flux divergence can be computed directly with a central difference scheme since numerical flux is smooth enough.
  Similar to WENO5 reconstruction, WENO5 interpolation calculates upwind $u_{i+1/2}^{\mathrm{left}}$ by weighting three interpolations of the third order.
  \begin{eqnarray}
    u_{i+1/2}^{\mathrm{left}} &=  \omega_1 C_1(u_{i-2},u_{i-1},u_{i})
                 +\omega_2 C_2(u_{i-1},u_{i},u_{i+1})
                 +\omega_3 C_3(u_{i},u_{i+1},u_{i+2})
  \end{eqnarray}
  where $C_j, j=1\cdots 3$ is the interpolation of the substencil defined by $u_{i+j-3}$, $u_{i+j-2}$ and $u_{i+j-1}$, while $\omega_j$ is the nonlinear weight which is samely defined as the weight in WENO5 reconstruction.
  More details about WENO5 interpolation is given in Deng's work\cite{deng_developing_2000}.

  We think that this process (WENO interpolation, solving Riemann flux and applying central difference) should be still valid in unstructured meshes.
  With some kind of WENO interpolation, we first compute the numerical flux on each edge middle point.
  And then the flux divergence on each vertex can be computed directly without numerical dissipation since numerical flux on edges is smooth enough, which should not be difficult.
  The first task is to effectively apply WENO interpolation on a stencil assmebled with connected edges.

  The WENO5 interpolation  in WCSN-E5 scheme is designed for uniform one-dimensional grid.
  Initial tests showed feasibility of this idea\cite{kunqu_2020,2021arXiv210212933Z}.
  But the accuracy is not satisfying since the stencil in Fig. \ref{fig:assemble-stencil} is nonuniform.
  In order to keep the high accuracy, we tried to construct a smooth curve across all the 6 nodes and the middle point
  and then perform WENO5 interpolation based on the curvelength.
  Unfortuninately, the rate of convergence is still the first order.

\begin{figure}[hbtp]
    \centering
    \includegraphics[width=0.4\textwidth]{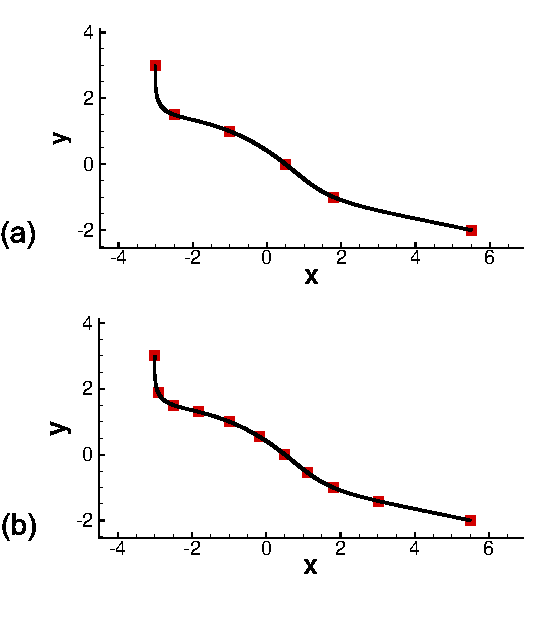}
    \caption{The interpolation for the ``inserting'' case: (a) a baseline curve and the interpolation nodes; (b) the refined curve.}
    \label{nuweno_refine_big}
  \end{figure}

\begin{figure}[hbtp]
    \centering
    \includegraphics[width=0.4\textwidth]{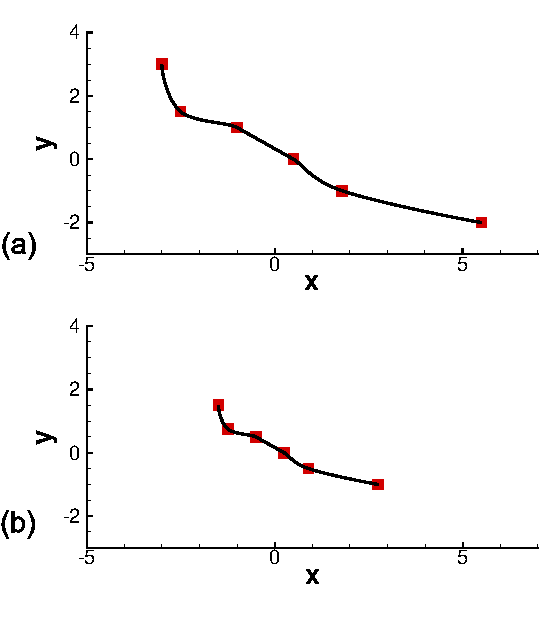}
    \caption{The interpolation for the ``shrunk'' case: (a) a baseline curve and the interpolation nodes; (b) the refined curve.}
    \label{nuweno_scale_big}
  \end{figure}

  After some analysis, we belive that the style of mesh refinement plays an important role here.
  Two kinds of refinements are involved in CFD.
  In the first kind, as shown in Fig.\ref{nuweno_refine_big}, the same grid curve is discretized with more and more nodes.
  This is just the style of refining a structured grid.
  In the second kind of refinement shown in Fig.\ref{nuweno_scale_big}, the mesh (and the edges in the mesh) shrinks agian and agian, with the interpolation nodes move closer and closer to increase resolution of the mesh.
  This style is similar to the refinement of unstructured meshes in which elements become smaller and smaller.

  Consider the $p$-th order Lagrangian interpolation in the $\xi$ domain on a curve 
  The leading term of the truncated error is
  \begin{equation}
    \epsilon = \frac{1}{p!}\frac{\mathrm{d}^p u}{\mathrm{d} \xi^p} (\Delta \xi)^p
  \end{equation}
  Along this 2D curve, given the transformation $x(\xi)$ and $y(\xi)$, we can transform this term into $x-y$ plane.
  The transformed leading error must contain a term
  \begin{equation} \label{eq:err_term}
  \frac{1}{p!}\frac{\partial^p u}{\partial x^p} \frac{\mathrm{d}^p x}{\mathrm{d} \xi^p}  (\Delta \xi)^p
  \end{equation}

  If the curve is refined with factor $N$ by means of the first kind of refinement, we can set that $\Delta \hat{\xi} = \Delta \xi /N $ on the refined curve.
  Thus
  \begin{equation} \label{eq:err_term_1}
    \frac{1}{p!}\frac{\partial^p u}{\partial x^p} \frac{\mathrm{d}^p x}{\mathrm{d} \xi^p}  (\Delta \xi)^p
    = \frac{1}{ N^p}\:\frac{1}{p!}\frac{\partial^p u}{\partial x^p} \frac{\mathrm{d}^p x}{\mathrm{d} \xi^p}  (\Delta \xi)^p
  \end{equation}
  which means the error decreases with a factor $1/N^p$.
  Thus it is of the $p$-th order.

  However, refining an unstructred mesh just likes shrinking the patterns in the mesh.
  When shrinking the length scale from $\Delta x$ to $\Delta \hat{x}$, we have the new transformation
  \begin{equation}\label{equ:scaledtrans}
    \hat{x}(\xi) = \dfrac{1}{N}x(\xi) \qquad
    \hat{y}(\xi) = \dfrac{1}{N}y(\xi)
  \end{equation}
  where the scaling factor $N>1$.
  From this, the relationship between derivatives of the two curves
  \begin{equation} \label{eq:scaled_derivative}
    \dfrac{\mathrm{d}^m \hat{x}}{\mathrm{d} \xi^m} = \dfrac{1}{N} \dfrac{\mathrm{d}^m {x}}{\mathrm{d} \xi^m} \qquad
    \dfrac{\mathrm{d}^m \hat{y}}{\mathrm{d} \xi^m} = \dfrac{1}{N} \dfrac{\mathrm{d}^m {y}}{\mathrm{d} \xi^m}
  \end{equation}
  The variation in the shirinked mesh must be
  \begin{equation} \label{eq:scaled_dx}
    \Delta \hat{x} = \dfrac{1}{N}\Delta x \qquad
    \Delta \hat{y} = \dfrac{1}{N}\Delta y
  \end{equation}
  It should be pointed out that, different from the first kind refinement, we have $\Delta \xi = \Delta \hat{\xi}$ in this case.
  Insert Eq. (\ref{eq:scaled_derivative}) and (\ref{eq:scaled_dx}) into Eq. (\ref{eq:err_term}),
  the term can be converted as
  $$
  \frac{1}{p!} \frac{\partial^p u}{\partial x^p}
  \frac{\mathrm{d}^p x}{\mathrm{d} \hat{\xi}^p}  (\Delta \hat{\xi})^p
  = \frac{1}{N}  \frac{1}{p!} \frac{\partial^p u}{\partial x^p}
  \frac{\mathrm{d}^p {x}}{\mathrm{d} \xi^p} (\Delta \xi)^p
  $$
  This term decreases with $1/N$, which means the rate of  convergence is the first order.
  It is the result of multi-dimensional essence of the curve stencil.
  Thus, when refining an unstructured mesh, the rate of convergence of interpolation based on the generalized coordinate $\xi$ is only first order no matter how large of the degree of the polynomial.

  \subsection{Two-Dimensional Hermite Interpolating on A Curvilinear Stencil}

  According to the analysis above, multidimensional interpolation  is neccesary in order to keep high order accuracy.
  %
  It is well known that multidimensional interpolation need much more supporting nodes than one-dimensional interpolation in case of the same degree.
  That is why the stencil of a multidimensional Lagarangin interpolation in FR or a multidimensional reconstruction in FVM always spreads many nodes or cells.
  We do not prefer to using such a large stencil since we still try to apply some kind one-dimensional computing to reduce computing cost and complexity.
  In order to introduce more supporting data for the multidimensional interpolation on the curilinear stencils assembled in unstructured meshes, we try to use more DOFs on each node, but not more nodes.

  \begin{figure}
    \centering
    \includegraphics[width=0.45\textwidth]{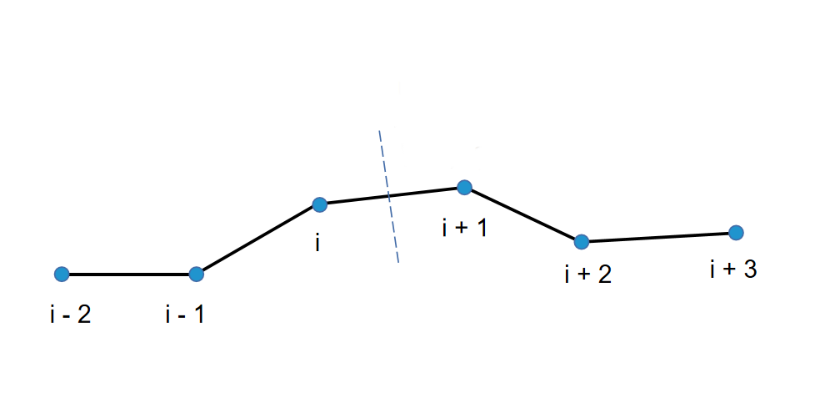}
    \caption{A schematic view of a stencil assembled in unstructured meshes.}
    \label{fig:sixnodes_stencil}
  \end{figure}
  In Fig. \ref{fig:sixnodes_stencil}, a stencil made of 5 edges and 6 nodes is shown.
  On the left and right sides of the middle point of the central edge betweeen Node $i$ and $i+1$, we try to interpolate the characteristic variables $\mathbf{c}_{i+1/2,L}$ and $\mathbf{c}_{i+1/2,R}$
  from which the numerical flux will be calculated along the direction of the central edge.
  Taking $\mathbf{c}_{i+1/2,L}$ for example, Nodes $i-2$, $i-1$, $\cdots$, $i+2$ are used.
  For the full stencil of 5 nodes, there are $15$ supporting data which satisfy the least number of conditions of 2D polynomial of the 4th degree.
  And these supporting datas give Equ.(\ref{eq:hintp_5}) to determine a polynomial
  \begin{equation}  \label{eq:hintp_5}
    \left\{
      \begin{array}{lcl}
     \mathbf{p}^T  (x_j,y_j) \:\mathbf{k} &=&  c_j  \\
     \mathbf{p}^T_x(x_j,y_j) \:\mathbf{k} &=&  c_{j,x}\\
     \mathbf{p}^T_y(x_j,y_j) \:\mathbf{k} &=&  c_{j,y}
      \end{array}
    \right.
    \qquad j=i-2, \cdots, i+2
  \end{equation}
  where
  \begin{eqnarray}
  \mathbf{p}(x,y)   &=& \left[1,x,y,\frac{x^2}{2},xy,\frac{y^2}{2},
                             \frac{x^3}{6},\frac{x^2y}{2},\frac{xy^2}{2},\frac{y^3}{6}, \right.  
                    \left. \frac{x^4}{24},\frac{x^3y}{6},\frac{x^2y^2}{4},\frac{xy^3}{6},\frac{y^4}{24}  \right]^T \\
  \mathbf{p}_x(x,y) &=& \left[0,1,0,
          x,y,0,
           \frac{{x}^{2}}{2}, xy,\frac{{y}^{2}}{2},  0, \right.  
                     \left. \frac{{{x}^{3}}}{6},\frac{{{x}^{2}}y}{2},\frac{x\,{{y}^{2}}}{2},\frac{{{y}^{3}}}{6},0 \right]^T\\
  \mathbf{p}_y(x,y) &=& \left[0,0,1,
                              0,x,y,
                              0,\frac{{x}^{2}}{2}, xy,\frac{{y}^{2}}{2}, \right.
                    \left.0, \frac{{{x}^{3}}}{6},\frac{{{x}^{2}}y}{2},\frac{x\,{{y}^{2}}}{2},\frac{{{y}^{3}}}{6} \right]^T
  \end{eqnarray}
  $c_j$ is any component of $\mathbf{c}_j$ and vector $\mathbf{k}$ contains all the coeffcients of the complete 4th degree 2D polynomial corresponding to the component.
  Collecting the equations on all 5 nodes, we can obtain a linear system of $15\times 15$.
  \begin{equation}
    \mathbf{Pk}=\mathbf{r}
  \end{equation}
  where
  \begin{equation}
    \mathbf{P} =
   \left[\begin{array}{c}
     \mathbf{p}^T  (x_{i-2},y_{i-2}) \\
     \mathbf{p}^T_x(x_{i-2},y_{i-2}) \\
     \mathbf{p}^T_y(x_{i-2},y_{i-2}) \\
     \vdots \\
     \mathbf{p}^T  (x_{i+2},y_{i+2}) \\
     \mathbf{p}^T_x(x_{i+2},y_{i+2}) \\
     \mathbf{p}^T_y(x_{i+2},y_{i+2})
   \end{array} \right]
   \;
   \mathbf{k} =
   \left[\begin{array}{l}
     k_1 \\
     k_2 \\
     \;\vdots \\
     k_{15}
    \end{array} \right]
   \;
   \mathbf{r} =
   \left[\begin{array}{l}
     c_{i-2} \\ c_{i-2,x} \\ c_{i-2,y} \\
     \;\;\vdots \\
     c_{i+2} \\ c_{i+2,x} \\ c_{i+2,y}
    \end{array} \right]
  \end{equation}
  Because matrix $\mathbf{P}$ might be ill-condistioned, we solve the linear system by means of SVD
  \begin{equation}
    \mathbf{P}  =  \mathbf{U}   \mathbf{\Sigma}  \mathbf{V}^T
    \quad \Longrightarrow \quad
    \mathbf{P}^{\dag}  = \mathbf{V} \mathbf{\Sigma} \mathbf{U}^T
    \quad \Longrightarrow \quad
    \mathbf{k} = \mathbf{P}^{\dag}   \mathbf{r}
  \end{equation}
  By ignoring some tiny singular values, we can solve $\mathbf{k}$ to obtain the polynomial for interpolation with high order accuracy even for ill-condistioned problems.
  We select singular values through trial and error.
  After SVD decomposition, we initially keep all non-zero singular values.
  If the absolute value of any component of $\mathbf{k}$ is not less than $1$, the smallest singular value should be dropped.
  And then try again for the left singular values until the absolute value of each component of $\mathbf{k}$ is less than $1$.

  If the origin is set at the middle point, the interpolated value is just $k_1$ which can be expressed as
  \begin{equation}
    k_1 = \mathbf{P}^{\dag}_{(1,:)} \mathbf{r}
  \end{equation}
  where $\mathbf{P}^{\dag}_{(1,:)} $ is the first row-vector of $\mathbf{P}^{\dag}$.
  Thus we can just store $\mathbf{P}^{\dag}_{(1,:)} $.

\begin{figure}
    \centering
    \includegraphics[width=0.4\textwidth]{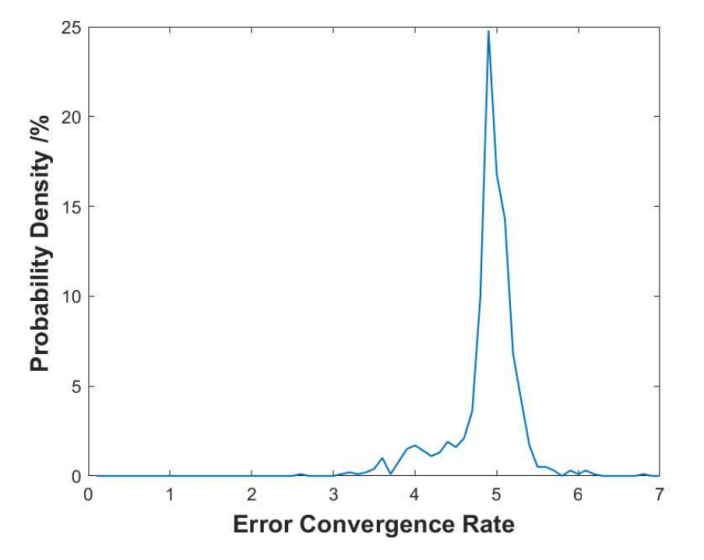}
    \caption{Probability Distribution of Convergence Order of Five-node Interpolation Error.}
    \label{fig:montocarol-5}
  \end{figure}

  Momte Carol test was performed to validate the accuracy of this Hermite interpolation.
  Stenciles of five-node with different length of edges and angles of deflection were generated randomly according to Gaussian distribution.
  The variance of edge length $\sigma_L$  is set as $15\%$ of the length of the middle edge.
  And for the angles of deflection between each two connected edges, the mean value $\mu_{\theta}=0$ and variance  $\sigma_{\theta}=\pi/12$.
  At the same time, different orientations of the stencil in the coordinate system were also considered.
  All the generated stencils were used to evaluate the order of accuracy.
  Fig.\ref{fig:montocarol-5} showed that this interpolation is of fifth order accuracy.

  \subsection{Hermite Finite Diference Method in 2D Unstructured Meshes Base on Hermite Interpolation}

  \subsubsection*{Governing Equation}

  In the above Hermite interpolation, gradient are neccesary.
  In order to avoid approximating gradient from some 2D stencils, we choose to solve conservative variable $\mathbf{w}$ and its gradient at the same time, just like Taylor DG\cite{luo_discontinuous_2008,luo_discontinuous_2009} or Hermite FD/FV\cite{zhu_new_2018,qiu_hermite_2004}.
  Thus, following equations are solve in our work
  \begin{equation}\label{eq:HWENO}
  \left\{
  \begin{array}{lll}
  \dfrac{\partial \mathbf{w}}{\partial t}   &+ \nabla \cdot \mathbf{f}  (\mathbf{w}) &= 0\\
  \\
  \dfrac{\partial \mathbf{w}_x}{\partial t} &+ \nabla \cdot \mathbf{f}_x(\mathbf{w}) &= 0\\
  \\
  \dfrac{\partial \mathbf{w}_y}{\partial t} &+ \nabla \cdot \mathbf{f}_y(\mathbf{w}) &= 0\\
  \end{array}
  \right.
  \end{equation}
  where $\mathbf{w}$ is the conservative variables, $\mathbf{w}_x$ and $\mathbf{w}_y$ are the first order derivatives (gradient) of $\mathbf{w}$.
  They are all unknowns and simutaniously evolved in time domain.
  In the course of interpolation, $\mathbf{w}$, $\mathbf{w}_x$ and $\mathbf{w}_y$ on each vertex/cell are converted to characteristic variables and used to perform Hermite interpolation.
  By this way, the stencil can be kept compact.

  \subsubsection*{Capturing Discontinuities}
  In the classical one-dimensional WENO5 interpolation, the nonlinear weights approach to the optimal values in smooth regions, generating the 5th order interpolation.
  However, due to the arbitrary distribution of nodes, it is difficult to obtain such optimized weights analytically in our case.
  For the sake of simplicity, the full 5th order interpolation is used if the solution is smooth, otherwise the nonlinear weighted interpolation is used.

  Here the three weights of the one-dimensional  uniform  WENO5 interpolation are first computed.
  The weight of the first sub-stencil, $\omega_1$, is used to determine whether the stecil is in discontinuous regions.
  Our numerical experiments show that $\omega_1$ ranges from $0.01$ to $0.79$ when the stencil is located in smooth regions.
  Then the task is to obtain a non-oscillation low order interpolation if the stencil is in discontinuous regions.
  In our numerical tests, two-dimensional 3-node Hermite interpolation makes the evolution of $ \mathbf{w}_x$ and  $ \mathbf{w}_y$ unstable.
  Referring to Zhu's work\cite{zhu_new_2018}, we chose to ignore the derivative conditions when interpolating near the discontinuities.
  This means that Lagrangian interpolation but not Hermite interpolation is used.
  The assumed polynomial is linear in two-dimensional space since we have only three function values on the three nodes of each sub-stencil.
  This results in the 2nd order accuracy near discontinutities.
  Thus the interpolation can be summerized as
  \begin{eqnarray}\label{eq:2d_hweno_intp}
    C^{\mathrm{left}}_{i+1/2} = \left\{
     \begin{array}{cl}
      C^{(5)}   &  \quad \mathrm{if\:} w_1 \in (0.01,\: 0.79) \\
       & \\
       \sum\limits_{j=1}^{3} \omega_j C^{(2)}_{j}   &  \quad\mathrm{else}
     \end{array}
    \right.
  \end{eqnarray}
  where $C^{(2)}_{j}$ is the value of a Lagrangian interpolation of the 2nd order , while $C^{(5)}$ is the value obtained with the 5th order Hermite interpolation.
  %

  %
  %
  For example, a three-node interpolation can be determined by
  \begin{eqnarray}\label{eq:3nodes_intp}
    \left[ \begin{array} {l}
      \mathbf{p}^T_3(x_{i-2},y_{i-2}) \\
      \mathbf{p}^T_3(x_{i-1},y_{i-1}) \\
      \mathbf{p}^T_3(x_{i}\;\;\:\:,y_{i}\;\;\:\:)
    \end{array}  \right]
    \left[ \begin{array} {c}
      k_1 \\ k_2 \\ \vdots \\ k_6
    \end{array}  \right]=
    \left[ \begin{array} {l}
      C_{i-2} \\  C_{i-1} \\  C_{i}
    \end{array}  \right]
  \end{eqnarray}
  where
  \begin{equation}
    \mathbf{p}^T_3 = \left[1,\: x,\: y, \: \frac{x^2}{2},\: xy,\: \frac{y^2}{2} \right]
  \end{equation}
  Again, the LHS matrix of $3 \times 6$ can be also decomposed with SVD.
  And thanks for SVD, we can obtain interpolation coefficients even the linear system is underdetremined.
  If the three nodes are collinear, it is degraded into a one-dimensional second degree polynomial interpolation which is the third-order accuracy.
  Otherwise, we obtain a two-dimensional linear interpolation with second-order accuracy.
  Algoritm \ref{alg:algorithm1} summarizes the procedure of this two-dimensional WENO interpolation along a five-node stencil.

  \begin{algorithm}[t]
    \caption{Hermite WENO5 interpolation on an edge}
    \label{alg:algorithm1}
    \LinesNumbered
    \KwIn{$[\mathbf{x}_{i-2}, \cdots \mathbf{x}_{i+2}]$,
          $[\mathbf{w}_{i-2}, \cdots \mathbf{w}_{i+2}]$,
          $[k_1,\cdots,k_{15}]$,
          and $\mathbf{n}$ }
    \KwOut{$\mathbf{w}^{\mathrm{left}}_{i+1/2}$ on the middle point of the edge.}
    \BlankLine
    Compute Roe average $\tilde{\mathbf{w}}$ from $\mathbf{w}_i$ and  $\mathbf{w}_{i+1}$   \;

    Compute $\mathbf{L}$ and  $\mathbf{R}$ from  $\tilde{\mathbf{w}}$ and $\mathbf{n}$;

    Obtain characteristic variables $\mathbf{c}_k = \mathbf{L} \:\mathbf{w}_k$ where $k=i-2,\cdots, i+2$;

    \ForEach{$k$-th component of characteristic variables}{
      Compute the nonlinear weights $w_1$, $w_2$ and $w_3$ by means of the uniform WENO5 interpolation;

      \eIf{$w_1 \in (0.01,0.79)$}
      {
        Solve three interpolations like Eq.(\ref{eq:3nodes_intp}) to obtain $C^{(2)}_{j}$;

        Compute $C^{\mathrm{left}}_{i+1/2,j}$ by means of nonlinear weighted summation;
      }
      {
        Compute $C^{\mathrm{left}}_{i+1/2,j}$ as $C^{(5)}$ with $[k_1, \cdots, k_{15}]$;
      }
    }

    \Return  $\mathbf{w}^{\mathrm{left}}_{i+1/2}= \mathbf{R}\: \mathbf{c}^{\mathrm{left}}_{i+1/2}$;

  \end{algorithm}

  With this interpolation, we can obtain $\mathbf{w}^{\mathrm{left}}_{i+1/2}$ and  $\mathbf{w}^{\mathrm{right}}_{i+1/2}$ at the middle point of the edge.
  Then the numerical flux $\mathbf{h}_{i+1/2}$ along the edge at the middle point $i+1/2$ can be computed with a Riemann flux solver.
  The numerical flux  $\mathbf{h}_{i+1/2}$ is the project of $[\mathbf{f},  \mathbf{g}]_{i+1/2}^T$ along the direction of the edge.
  At the same time, full vectors $[\mathbf{f},  \mathbf{g}]^T$   at each vertex can be analytically computed on each vertex since its $\mathbf{w}$ is known.

  \subsubsection*{Computing Divergence with Least-Square Method}
  \begin{figure}
    \centering
    \includegraphics[width=0.4\textwidth]{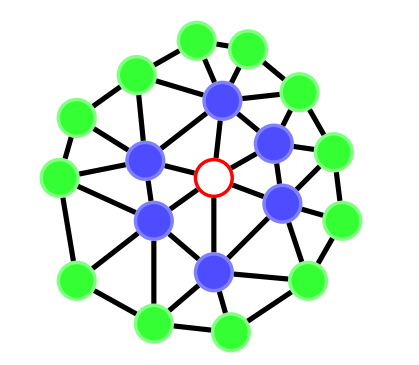}
    \caption{The double level stencil used to compute flux divergence on the centeal node (noted as a red circle). The blue nodes and green nodes are the first level and the second level neighbors, respectively. Each node holds $f$ and $g$, while the edges have projection of $[f,g]^T$.}
    \label{fig:div-stecil}
  \end{figure}
In order to compute the divergence, a double level stencil around a node can be constructed, as shown in Fig.\ref{fig:div-stecil}.
It contains two level of nodes and all the edges between the nodes.
The flux data on these nodes and edges can be used to derive two polynomials of the flux fields from which the divergences in Eq. (\ref{eq:HWENO}) can be computed.

  For a couple of components $f$ and $g$ of $[\mathbf{f},\mathbf{g}]^T$, we approximate them as two-dimensional polynomials
  \begin{equation}
    \left \{
      \begin{aligned}
        f(x,y) &= \mathbf{p}^T(x,y) \mathbf{a} \\
        g(x,y) &= \mathbf{p}^T(x,y) \mathbf{b}
      \end{aligned}
    \right.
  \end{equation}
  where $\mathbf{p}$ contains terms up to specified degree.
  For a given vertex, we can set a linear system from flux datas on vertices and edges  arround $Q$.

  \begin{equation}
    \label{eq:flux_polynomials}
    \left[ \begin{array}{c|c}
      \mathbf{p}^T(x_1,y_1)  & \mathbf{O} \\
      \mathbf{O} & \mathbf{p}^T(x_1,y_1) \\
      \vdots & \vdots \\
      \mathbf{p}^T(x_n,y_n) &  \mathbf{O}\\
      \mathbf{O} & \mathbf{p}^T(x_n,y_n) \\
      \vdots & \vdots\\
      \mathbf{p}^T(x_N,y_N) &  \mathbf{O} \\
      \mathbf{O} & \mathbf{p}^T(x_N,y_N) \\
      \hline \\
      \mathbf{p}^T(\bar{x}_{1},\bar{y}_{1}) n_{x,1}  &
      \mathbf{p}^T(\bar{x}_{1},\bar{y}_{1}) n_{y,1} \\
      \vdots &\vdots \\
      \mathbf{p}^T(\bar{x}_{m},\bar{y}_{m}) n_{x,m}  &
      \mathbf{p}^T(\bar{x}_{m},\bar{y}_{m}) n_{y,m} \\
      \vdots &\vdots \\
      \mathbf{p}^T(\bar{x}_{M},\bar{y}_{M}) n_{x,M} &
      \mathbf{p}^T(\bar{x}_{M},\bar{y}_{M}) n_{y,M}
    \end{array}\right]
    \left[
      \begin{array}{c}
        \mathbf{a} \\
        \mathbf{b}
      \end{array}
    \right] =
    \left[
      \begin{array}{l}
         f_1 \\
         g_1 \\
         \;\vdots \\
         f_n \\
         g_n \\
         \;\vdots \\
         f_N \\
         g_N\\
         \hline \\
         h_{1} \\
         \;\vdots \\
         h_{m}\\
         \;\vdots \\
         h_{M}
      \end{array}
    \right]
  \end{equation}
  In the equations above, the first $2N$ equations come from $N$ vertices, while the left $M$ equations come from $M$ edges.
  $(x_n,y_n)$ is the location of the $n$-th vertex.
  And  $(\bar{x}_m,\bar{y}_m)$ is the location of the middle point of  the $m$-th edge.
  $[n_{x,m},n_{y,m}]$ is the unit direction of the $m$-th edge.

  The verties and edges in the stencil make Eq. (\ref{eq:flux_polynomials})  an over-determined system.
  And it can be solved with SVD.
  \begin{equation}
    \left[
      \begin{array}{c}
        \mathbf{a} \\
        \mathbf{b}
      \end{array}
    \right] = \mathrm{LHS}^{\dag} \; \mathrm{RHS}
  \end{equation}
  Agian, since  we can directly express such divergences with some components of $\mathbf{a}$ and $\mathbf{b}$ by defining the central node  as the origion point.
  \begin{equation}
    \left \{
      \begin{aligned}
        \nabla \cdot \mathbf{f}(0,0) &= \left[f_x +  g_y \right]_{(0,0)}
        = a_2 + b_3 \\
        &=  \left[(\mathrm{LHS}^{\dag})_{(a_2,:)} +  (\mathrm{LHS}^{\dag})_{(b_3,:)}\right] \: \mathrm{RHS}\\
        & \\
        \nabla \cdot \mathbf{f}_x(0,0) &= \left[ \frac{\partial f_x}{\partial x} +  \frac{\partial g_x}{\partial y} \right ]_{(0,0)} = a_4 + b_5\\
        &=  \left[(\mathrm{LHS}^{\dag})_{(a_4,:)} +  (\mathrm{LHS}^{\dag})_{(b_5,:)}\right] \: \mathrm{RHS}\\
        & \\
        \nabla \cdot \mathbf{f}_y(0,0) &= \left[ \frac{\partial f_y}{\partial x} +  \frac{\partial g_y}{\partial y}\right]_{(0,0)}
        = a_5 + b_6\\
        &=  \left[(\mathrm{LHS}^{\dag})_{(a_5,:)} +  (\mathrm{LHS}^{\dag})_{(b_6,:)} \right] \:\mathrm{RHS}
      \end{aligned}
    \right.
  \end{equation}
  Thus only three row vectors,
  \begin{equation}
    \begin{aligned}
    \mathbf{d}_0 &= (\mathrm{LHS}^{\dag})_{(a_2,:)} +  (\mathrm{LHS}^{\dag})_{(b_3,:)} \\
    \mathbf{d}_x &= (\mathrm{LHS}^{\dag})_{(a_4,:)} +  (\mathrm{LHS}^{\dag})_{(b_5,:)} \\
    \mathbf{d}_y &= (\mathrm{LHS}^{\dag})_{(a_5,:)} +  (\mathrm{LHS}^{\dag})_{(b_6,:)}
    \end{aligned}
  \end{equation}
  have to be precomputed and stored.

  \subsubsection*{Computing Viscous Flux}
  For viscous problems, we can compute components of viscous flux on each vertex directly
  \begin{equation}
    \mathbf{f}_{\mathrm{vis},i} = \mathbf{f}_{\mathrm{vis}} \left(
        \mathbf{w}_i,
        \left. \frac{\partial \mathbf{w}}{\partial x}\right|_{i} ,
        \left. \frac{\partial \mathbf{w}}{\partial y}\right|_{i}
        \right)
  \end{equation}
  After that, the divergence of $\mathbf{f}_{\mathrm{vis}}$ on a vertex can be also computed from the viscous flux data on the surround vertices by means of least-square method.
  Using the vertices in the stencil of the divergence of the inviscid flux, the forth order accuracy can be approached in two-dimensional cases.

  \vspace{100pt}\par

  With the divergences obtained, the spatial discretization is accomplished.
  Here the three-stage TVD Runger-Kutta is applied as the time scheme.

  \section{Numerical Tests}
  \subsection{Inviscid Isentropic Vortex}
  The two-dimensional moving isentropic vortex problem was adopted to evaluate the accuracy of our method.
  The initial field is defined by:
  \begin{equation}
    \left \{
    \begin{aligned}
      \rho &= \left[1-\frac{(\gamma-1)\beta^2}{8 \gamma \pi^2} \exp(1-r^2)\right]^{\frac{1}{\gamma-1}} \\
      (u,v) &= (1,1) + \frac{\beta}{2\pi} \exp\left(\frac{1-r^2}{2}\right) (-\bar{y},\bar{x}) \\
      p &= \rho ^ {\gamma}
    \end{aligned}
    \right.
  \end{equation}
  where $\beta =5$ and the specific heat ratio of the gas $\gamma=1.4$.
  This initial condition gives a vortex whose center is located at the origin.
  The domain is a square with $x\in [-10,+10]$ and $y\in [-10,+10]$.
  And the periodic boundary condition is imposed on its four edges.

  In order to test  the convergence order of the error with refined meshes, errors are measured on four different meshes.
  The coarsest mesh has $50$ unifom segments on each edge of the square domain.
  Delaunay triangulation algorithm is applied to generate an isentropic triangle mesh with 3957 nodes.
  Dividing the edges by $75$, $100$ and $150$ uniform segments respectively, another three isentropic triangle meshs can be obtained.
  Their node sizes are $8986$, $16163$ and $36408$.

\begin{figure}
    \centering
    \includegraphics[width=0.65\textwidth]{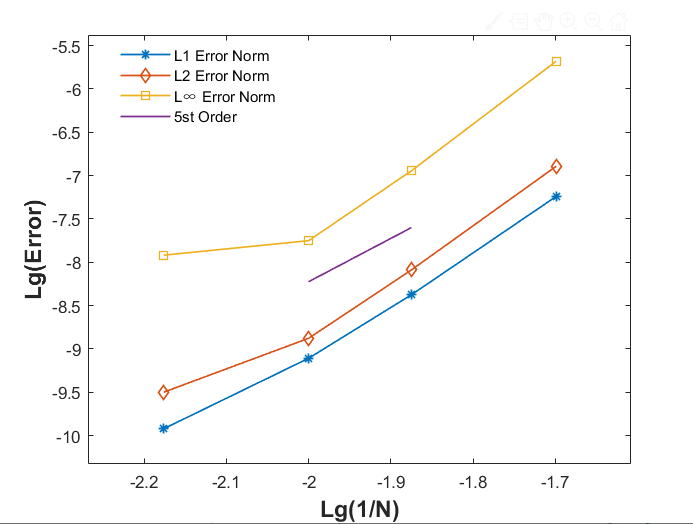}
    \caption{Norms of the density error \emph{vs} the length scale of the meshes. Here $N$ is the number of segments on each edge of the square domain.}
    \label{fig:isentropic-vortex}
  \end{figure}

  Fig. \ref{fig:isentropic-vortex} shows how the different norms of the density error decrease with the length scale of the meshes.
  We can see the order of accuracy is about fifth for this smooth problem.

  \subsection{Couette Flow}
  The two-dimensional compressible Couette flow refers to the flow between two parallel planes with different temperatures and constant relative movement.
  Becuase analytical solution exists, it is used to measure the order of accuracy for compressible viscous flows.

  Given two planes paralle to $x$ direction and the distance between them as $L=10$, as well as the speed of the top wall as $u_0$ with the bottom wall fixed, the exact solution is:
  \begin{equation}
     \left\{
      \begin{aligned}
        u &= \frac{2y}{L}u_0 \\
        T &= -\frac{1}{K} \left(\frac{2u^2_0}{L^2}\mu y^2 +C_1 y + C_0\right)\\
        \rho &= \frac{\gamma p}{(\gamma-1)C_p T}
      \end{aligned}
    \right.
  \end{equation}
  where
  \begin{equation}
    \left\{
     \begin{aligned}
       C_1 &= K\frac{T_{\mathrm{top}}-T_{\mathrm{bottom}}}{L}\\
       C_0 &= \frac{K(T_{\mathrm{top}}+T_{\mathrm{bottom}}) + \mu u^2_0}{2}\\
     \end{aligned}
   \right.
  \end{equation}

\begin{figure}
    \centering
    \includegraphics[width=0.45\textwidth]{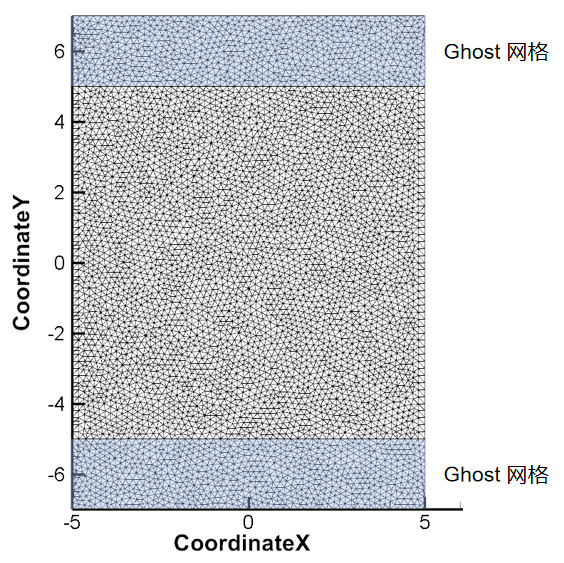}
    \caption{Domain and mesh of simulation of 2D Couette flow.}
    \label{fig:couette-mesh}
  \end{figure}

\begin{figure}
    \centering
    \includegraphics[width=0.45\textwidth]{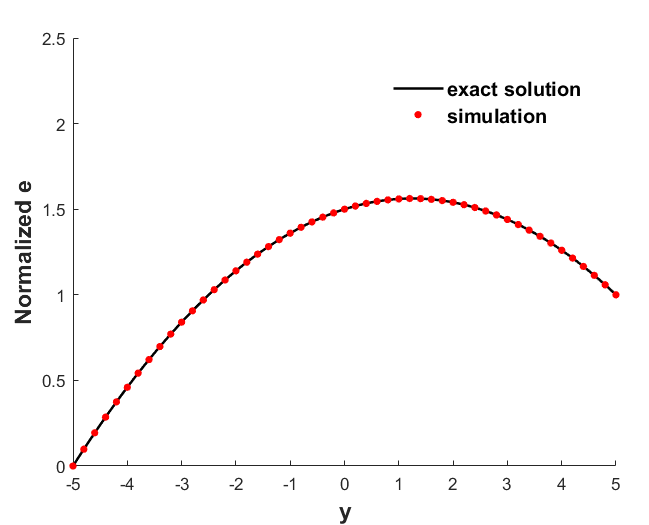}
    \caption{The nomalized temperature profile between the two planes of 2D Couette flow (computed with the mesh of $5556$ nodes).}
    \label{fig:couette-profile}
  \end{figure}

\begin{figure}
    \centering
    \includegraphics[width=0.5\textwidth]{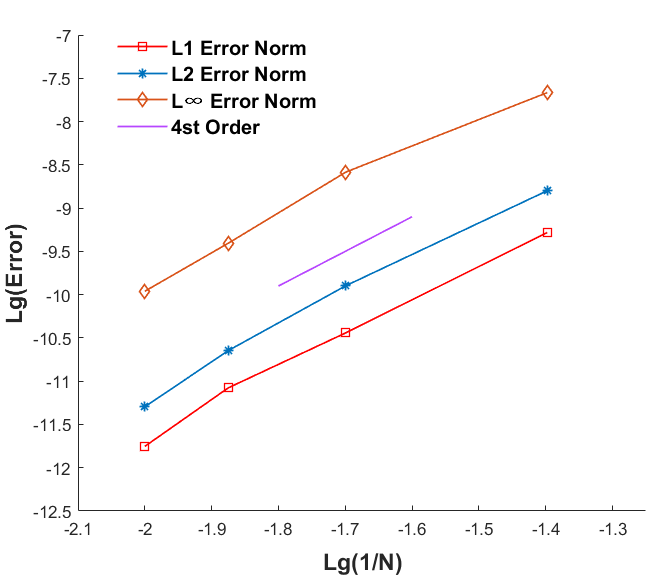}
    \caption{Mesh of simulation of 2D Couette flow.}
    \label{fig:couette-accuracy}
  \end{figure}

  The flow domain is defined as a square whose $x \in [-5,+5]$ and $y \in [-5,+5]$.
  Periodic boundary condition is imposed on the left and right edges.
  In order to impose the non-slip wall condition on the top and bottom walls, the domain is extended to define two ghost regions.
  Thus $y$ ranges from $-7$ to $+7$ (shown in Fig. \ref{fig:couette-mesh}).
  By unifomrly splitting each edge of the square domain into $35$, $70$, $105$ and $140$ segments and applying Delaunay triangulation algorithm, four isentropic triangle meshes with $1420$, $5556$, $12349$ and $22429$ nodes are generated for accuracy test.
  Simulations are successfully performed in the four meshes by setting $\gamma=1.4$, heat capacity $C_p=10$, heat conductivity $K=1/100$, vciscosity $\mu=1/10$.
  Fig. \ref{fig:couette-profile} shows the nomalized temperature profile between the two planes which agrees the exact solution perfectly.
  Fig. \ref{fig:couette-accuracy} shows the density error converges as the length scale of the meshes decreases with the forth order.
  %

  \subsection{Shock in 2D Space}
  The two-dimensional circular Sod problem which also has exact solution,  is simulated in order to test the capability to capture discontinuities of our method.
  The initial condition is defined as:
  \begin{equation}
    (\rho_0,u_0,v_0,p_0) = \left\{
      \begin{array}{cc}
        (1,\:0,\:0,\:1)  \qquad &  \mathrm{if}\: |r|<0.5 \\
        \left(\frac{1}{8},\:\frac{1}{4},\:0,\:\frac{1}{10}\right)  \qquad & \mathrm{else}
      \end{array}
    \right.
  \end{equation}
  where $\gamma=1.4$.

\begin{figure}
    \centering
    \includegraphics[width=0.35\textwidth]{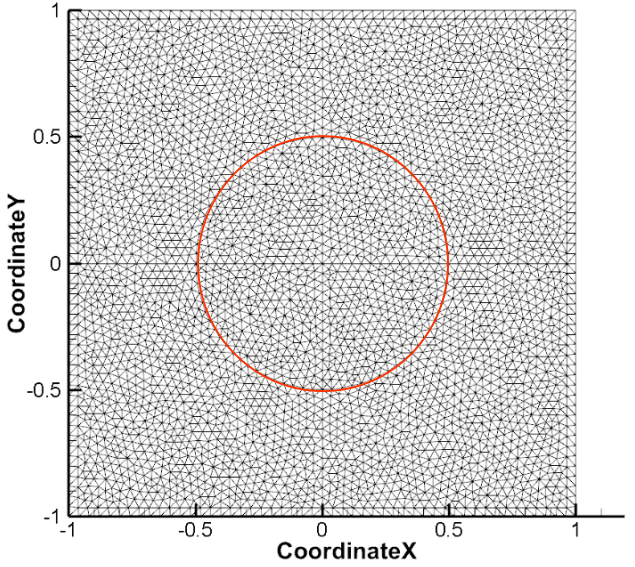}
    \caption{Mesh of 2D Sod case. The red circle is the boundary of the embedded region so that the initial field can be set perfectly.}
    \label{fig:sod2d-a-mesh}
  \end{figure}

  The meshes are still generated from unifomrly split four edges of a square of $2 \times 2$.
  But in order to perfectly set the initial field, a circular region of $r=0.5$ is embedded at the center of the square (Fig. \ref{fig:sod2d-a-mesh}).
  By unifomrly splitting each edge into $50$, $100$ and $150$ segments, we obtain four isentropic triganle meshes with $3970$, $15800$, $35458$ nodes respectively.

  \begin{figure}
    \centering
    \includegraphics[width=0.45\textwidth]{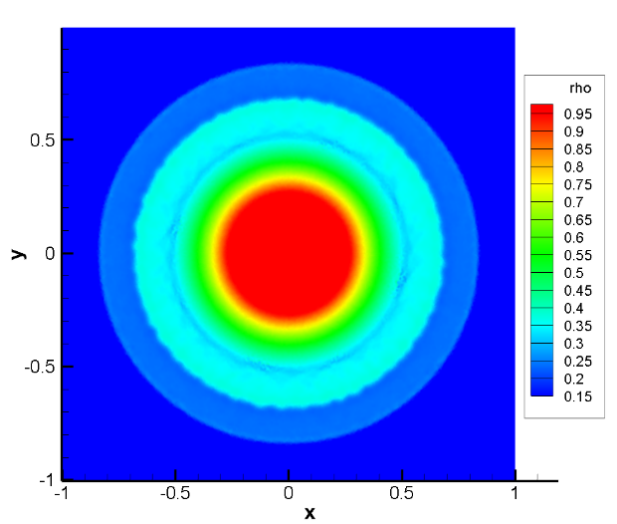}
    \includegraphics[width=0.45\textwidth]{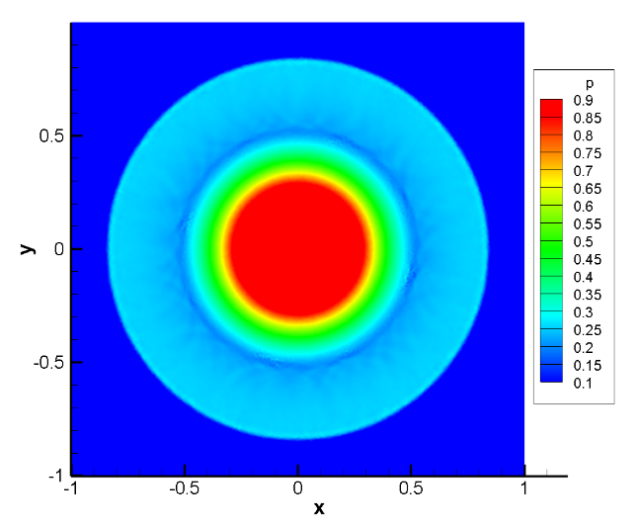}
    \caption{Contours of density(left) and pressure(right) computed at $t=0.2$ with mesh of 15800 nodes.}
    \label{fig:sod2d-a-contours}
  \end{figure}
\begin{figure}
    \centering
    \includegraphics[width=0.4\textwidth]{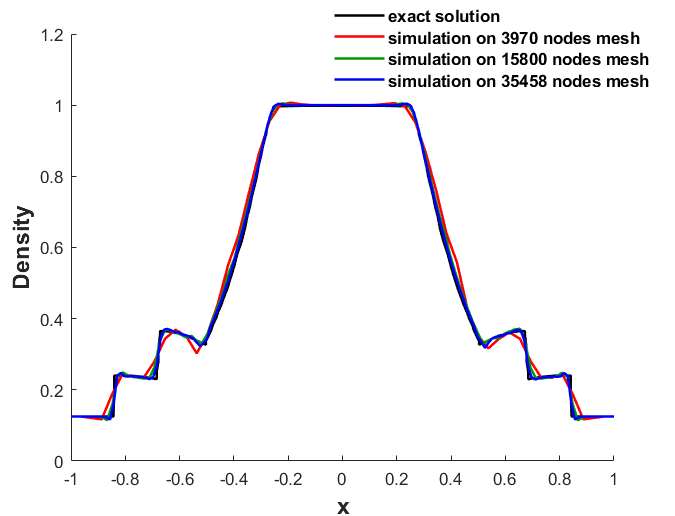}
    \includegraphics[width=0.4\textwidth]{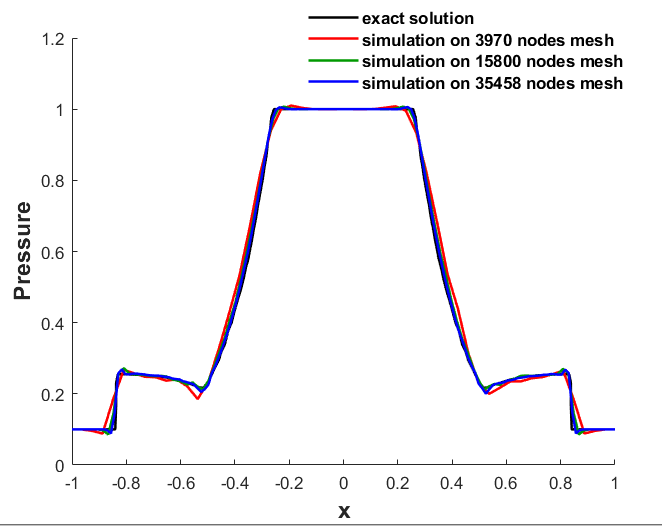}
    \caption{Profiles of density(top) and pressure(bottom) computed with the three meshes at $y=0$ and $t=0.2$, as well as the exact result presented as black profiles.}
    \label{fig:sod2d-a-profile}
  \end{figure}

  Fig. \ref{fig:sod2d-a-contours} presents the  contours  at $t= 0.2$ of the mesh with 15800 nodes.
  Fig. \ref{fig:sod2d-a-profile} presents the profiles of density and pressure computed in the three meshes with different resolution on $y=0$ at $t=0.2$, as well as the exact solution.
  It can seen that shocks and contact discontinuities are captured excellently.

  \section{Conclussion and Future Works}
  In this paper, in order to reduce the computing cost and complexity of capturing discontinuities in high order methods of unstructred meshes, the authors tried to extend the mature and high efficient one-dimensional WCNS FD scheme to unstructred meshes.

  A simple and economical WENO interpolation was proposed for two-dimensional unstructred triangle meshes.
  In smooth regions, to obtain the left and right values on the two sides of the middle point of each edge, the authors applied the two-dimensional Hermite interpolation along a curvilinear stencil which is assembled  with several connected edges in two-dimensional unstructred meshes.
  By means of SVD with selected singular values, the accuracy of this interpolation can be stably up to the fifth order on a stencil of only five nodes with not only the value but also the gradient on each node.
  Near discontinutities, by weighting three two-dimensional three-node Lagrangian interpolations, the second order accuracy can be obtained.
  The three three-node stecils are just the three substencils extracted from the original five-node stencil, which follows the same idea of WENO and WCNS schemes.
  And the weights also come from the one-dimensional WENO interpolation based on unifom grid.
  Because the interpolations are only applied along curvilinear stencils and gradient data are also used, the stencils are more compact and thus the interpolations are more economical.

  With the interpolated left and right values on each edge, the Riemann solver is applied to obtain a smooth inviscid flux.
  At the same time, analytical inviscid and viscous flux can be calculated directly on each vertex.
  From these flux data on each vertex and edge, the divergence of flux on each vertex can be approximated by means of the least-square based difference method.

  Monte Carol test on many random generated stencils showed the fifth order accuracy of this method for smooth fields.
  The isentropic vortex case and the Couette flow case show that this method is the fifth and forth order accracy for smooth inviscid and viscous problems respectively.
  At the same time, the results of the circular two-dimensional Sod problem show the validity of this method when capturing discontinuities in compressible flows.

  Because this paper consentrates on the spatial discretization, boundary conditions were not considered here.
  Thus most of the numerical tests are of periodic domains.
  At the same time, the least square based difference scheme adopted in this work is not conservative when computing the flux divergence.

  How to apply boundary conditions in this vertex-based method and how to modify this method to achieve conservation are new challenges.
  And extending this idea to three-dimensional problems is also an interesting topic.

  \bibliographystyle{unsrt}
  \bibliography{main}

\begin{thebibliography}{10}

\bibitem{huynh_flux_2007}
H.~T. Huynh.
\newblock A {Flux} {Reconstruction} {Approach} to {High}-{Order} {Schemes}
  {Including} {Discontinuous} {Galerkin} {Methods}.
\newblock In {\em 18th {AIAA} {Computational} {Fluid} {Dynamics} {Conference}},
  Miami, Florida, June 2007. American Institute of Aeronautics and
  Astronautics.

\bibitem{liu2006spectral}
Yen Liu, Marcel Vinokur, and Zhi~Jian Wang.
\newblock Spectral difference method for unstructured grids i: basic
  formulation.
\newblock {\em Journal of Computational Physics}, 216(2):780--801, 2006.

\bibitem{wang2007spectral}
Zhi~Jian Wang, Yen Liu, Georg May, and Antony Jameson.
\newblock Spectral difference method for unstructured grids ii: extension to
  the euler equations.
\newblock {\em Journal of Scientific Computing}, 32(1):45--71, 2007.

\bibitem{cockburn1989tvb}
Bernardo Cockburn and Chi-Wang Shu.
\newblock Tvb runge-kutta local projection discontinuous galerkin finite
  element method for conservation laws. ii. general framework.
\newblock {\em Mathematics of computation}, 52(186):411--435, 1989.

\bibitem{cockburn2001runge}
Bernardo Cockburn and Chi-Wang Shu.
\newblock Runge--kutta discontinuous galerkin methods for convection-dominated
  problems.
\newblock {\em Journal of scientific computing}, 16(3):173--261, 2001.

\bibitem{shu_efficient_1988}
Chi-Wang Shu and Stanley Osher.
\newblock Efficient implementation of essentially non-oscillatory
  shock-capturing schemes.
\newblock {\em Journal of Computational Physics}, 77(2):439--471, August 1988.

\bibitem{jiang_efficient_1996}
Guang-Shan Jiang and Chi-Wang Shu.
\newblock Efficient {Implementation} of {Weighted} {ENO} {Schemes}.
\newblock {\em Journal of Computational Physics}, 126(1):202--228, June 1996.

\bibitem{deng_developing_2000}
Xiaogang Deng and Hanxin Zhang.
\newblock Developing {High}-{Order} {Weighted} {Compact} {Nonlinear} {Schemes}.
\newblock {\em Journal of Computational Physics}, 165(1):22--44, November 2000.

\bibitem{kunqu_2020}
Kun Qu and Meiyuan Zhen.
\newblock A {High}-{Order} {Shock}-{Capturing} {Finite} {Difference} {Method}
  in {Unstructed} {Meshes}.
\newblock In {\em 11th {Chinese} {Computational} {Fluid} {Dynamics}
  {Conference}}, Shenzhen, Guangdong, December 2020. The Chinese Society of
  Theoretical and Applied Mechanics (CSTAM).

\bibitem{2021arXiv210212933Z}
Meiyuan {Zhen}, Kun {Qu}, and Jinsheng {Cai}.
\newblock {A Novel Finite Difference Method for Euler Equations in 2D
  Unstructured Meshes}.
\newblock {\em arXiv e-prints}, page arXiv:2102.12933, February 2021.

\bibitem{luo_discontinuous_2008}
Hong Luo, Joseph~D. Baum, and Rainald Löhner.
\newblock A discontinuous {Galerkin} method based on a {Taylor} basis for the
  compressible flows on arbitrary grids.
\newblock {\em Journal of Computational Physics}, 227(20):8875--8893, October
  2008.

\bibitem{luo_discontinuous_2009}
H.~Luo, J.D. Baum, and R.~Löhner.
\newblock A discontinuous {Galerkin} method using {Taylor} basis for computing
  shock waves on arbitrary grids.
\newblock In Klaus Hannemann and Friedrich Seiler, editors, {\em Shock
  {Waves}}, pages 1005--1010, Berlin, Heidelberg, 2009. Springer.

\bibitem{zhu_new_2018}
Jun Zhu and Jianxian Qiu.
\newblock New {Finite} {Volume} {Weighted} {Essentially} {Nonoscillatory}
  {Schemes} on {Triangular} {Meshes}.
\newblock {\em SIAM Journal on Scientific Computing}, 40(2):A903--A928, January
  2018.
\newblock Publisher: Society for Industrial and Applied Mathematics.

\bibitem{qiu_hermite_2004}
Jianxian Qiu and Chi-Wang Shu.
\newblock Hermite {WENO} schemes and their application as limiters for
  {Runge}–{Kutta} discontinuous {Galerkin} method: one-dimensional case.
\newblock {\em Journal of Computational Physics}, 193(1):115--135, January
  2004.

\end{thebibliography}

\end{document}